\documentclass[12pt,letterpaper]{amsart}
\usepackage{amssymb,latexsym, amsmath, amsxtra}
\usepackage[dvips]{graphics}

\theoremstyle{plain}
\newtheorem{theorem}{Theorem}[section]

\newtheorem{conjecture}[theorem]{Conjecture}
\newtheorem{question}[theorem]{Question}
\newtheorem{lemma}[theorem]{Lemma}

\theoremstyle{definition}

\theoremstyle{remark}

\numberwithin{equation}{section}
\numberwithin{theorem}{section}
\numberwithin{table}{section}
\numberwithin{figure}{section}

\newcommand{\C}{\mathbb C}

\newcommand{\R}{\mathbb R}

\renewcommand{\H}{\mathcal H}

\def\({\left(}
\def\){\right)}

\def \mm #1#2#3#4{\begin{pmatrix} #1 & #2 \cr #3 & #4 \end{pmatrix}}

\begin{document}
\title{$L$-functions and higher order modular forms}
\author{David W. Farmer and Sarah Zubairy}

\thanks{A portion of this work arose from an REU program at 
Bucknell University and the
American Institute of Mathematics.  Research supported by the
American Institute of Mathematics and the National Science Foundation}

\thispagestyle{empty}
\vspace{.5cm}
\begin{abstract}
It is believed that Dirichlet series with a functional equation and
Euler product of a particular form are associated to 
holomorphic newforms on a Hecke congruence group.  
We perform computer algebra experiments which
find that in certain cases one can associate a kind of
``higher order modular form'' to such Dirichlet series.
This suggests a possible approach to a proof of the conjecture.
\end{abstract}

\address{
{\parskip 0pt
American Institute of Mathematics\endgraf
farmer@aimath.org\endgraf
\null
Duke University\endgraf
}
  }

\maketitle

\section{Introduction}

We investigate the relationship between degree-2 $L$-functions and
modular forms.  
We find that degree-2 $L$-functions can be associated to functions on the
upper half-plane which have similar properties to ``second order
modular forms.''  Since it is conjectured that degree-2 $L$-functions
can be associated to modular forms, this looks like a step in the right direction.

We review some
classical results on modular forms and then describe 
the conjecture which motivates our work.
A good reference for this material is Iwaniec's book~\cite{Iw}.

Let 
$$
\Gamma_0(N)=\left\{\mm abcd : a,b,c,d \hbox{ are integers, } 
ad - bc = 1 \hbox{, and } c \equiv 0 \bmod N \right\}
$$
be the Hecke congruence group of level~$N$, and suppose 
$\chi$ is a character mod~$N$.  The group $\Gamma_0(N)$ acts
on functions $f:\mathcal H\to \C$ by $f\to f|\gamma$ where
\begin{equation}
f(z)\big| \mm abcd  = \chi(d)^{-1} (cz+d)^{-k} f\(\frac{az+b}{cz+d}\) .
\end{equation}
Here $\mathcal H=\{x+i y\in \C\ :\ y>0\}$ is the upper half of the
complex plane.
The vector space of \emph{cusp forms of weight $k$ and character $\chi$ for
$\Gamma_0(N)$}, denoted $S_k(\Gamma_0(N),\chi)$, is the set of 
holomorphic functions $f:\mathcal H\to \C$ which
satisfy $f|\gamma=f$ for all $\gamma\in \Gamma_0(N)$ and which
vanish at all cusps of~$\Gamma_0(N)$.  Since 
\begin{equation}
T := \mm 1101 \in \Gamma_0(N)
\end{equation}
we have $f(z)=f(z+1)$, so there is a Fourier expansion of the form
\begin{equation}\label{eqn:fourier}
f(z)=\sum_{n=1}^\infty a_n e^{2 \pi i n z} .
\end{equation}

In the case $\chi$ is the trivial character~$\chi_0$, the
newforms in $S_k(\Gamma_0(N),\chi_0)$ have a distinguished basis
of \emph{Hecke eigenforms} which satisfy
\begin{equation}\label{eqn:frickepm}
f|H_N=\pm f
\end{equation}
and
\begin{equation}\label{eqn:hecke}
f|T_p=a_p f
\end{equation}
for prime $p$.
Here
$$
H_N = \mm {}{-1}N{} 
$$
is the Fricke involution.  If $\ell$ is prime,
\begin{equation}\label{def:heckeop}
T_\ell = \chi(\ell) \mm \ell{}{}1 + \sum_{a=0}^{\ell-1} \mm 1a{}\ell,
\end{equation}
is the Hecke operator.  If $\ell | N$ then $\chi(\ell)=0$ and $T_\ell$ 
is known as the Atkin-Lehner operator~$U_\ell$.

We will now state our motivating conjecture, and then explain its relevance
to the theory of $L$-functions.

\begin{conjecture}\label{theconjecture} If  $f:\mathcal H\to \C$ is
analytic,  is periodic with period~1~(\ref{eqn:fourier}), and satisfies
the Fricke~(\ref{eqn:frickepm}) and Hecke~(\ref{eqn:hecke}) 
relations with $\chi=\chi_0$, then $f\in S_k(\Gamma_0(N),\chi_0)$.
\end{conjecture}
Thus, the invariance property $f|\gamma=f$, which leads to the
Fricke and Hecke relations, would actually follow from
them.

We will rephrase the conjecture in terms of $L$-functions.
Associated to a cusp form with Fourier expansion~(\ref{eqn:fourier})
is an $L$-function
\begin{equation}
L(s,f)=\sum_{n=1}^\infty \frac{a_n}{n^s} .
\end{equation}
Using the Mellin transform and its inverse, it can be shown that
the Fricke relation~(\ref{eqn:frickepm}) is 
equivalent to the functional equation
\begin{equation} \label{eqn:thefunctionalequation}
\xi_f(s):=\(\frac{2\pi}{\sqrt{N}}\)^{-s}\Gamma\(s\)L_f(s)=
\pm (-1)^{k/2}\xi_f(k-s).
\end{equation}
Also, the Hecke relations~(\ref{eqn:hecke}) are equivalent to
$L(s,f)$ having an Euler product of the form
\begin{equation} \label{eqn:theEulerproduct}
L(s,f)=\prod_p
\(1-a_p p^{-s} + \chi(p) p^{k-1-2s}\)^{-1},
\end{equation}
because both statements are equivalent to 
$a_{p^n m}=a_{p^n} a_m$ for $p\nmid m$ and
\begin{equation}
a_{p^{n+1}} = a_p a_{p^n} - \chi(p)p^{k-1} a_{p^{n-1}}.
\end{equation}
Thus, Conjecture~\ref{theconjecture} is equivalent to 
\begin{conjecture}\label{theconjectureL} If  a Dirichlet series
continues to an entire function of order one which is bounded in 
vertical strips and satisfies the functional 
equation~(\ref{eqn:thefunctionalequation}) 
and the Euler product~(\ref{eqn:theEulerproduct}) with $\chi=\chi_0$,
then the Dirichlet series equals $L(s,f)$ for some
$f\in S_k(\Gamma_0(N),\chi_0)$.
\end{conjecture}

This conjecture should be viewed as part of the Langlands'
program.  Note that one does not require functional equations for
twists of the $L$-function, as in Weil's converse theorem.
As a special case, the $L$-function of a rational elliptic curve
automatically has an Euler product of form~(\ref{eqn:theEulerproduct})
with $k=2$ and $\chi=\chi_0$,
so the modularity of a rational elliptic curve would 
follow from the analytic continuation and functional equation
for one $L$-function.  

Progress on the conjecture has been made only for small $N$, 
for the trivial character~\cite{CF}, and (appropriately modified)
for almost the
same cases for nontrivial character~\cite{Har}.
For $N\le 4$, Hecke's original converse theorem establishes
the conjecture.  This follows from the fact that the
group generated by $T$ and $H_N$ contains $\Gamma_0(N)$
exactly when $N\le 4$.  Note that this only uses the
functional equation, not the Euler product.
For larger~$N$, one must use the Euler product in a 
nontrivial way.  This possibility was introduced in~\cite{CF},
and examples were given for certain~$N\le 23$.

In this paper we specialize to the case~$N=13$, for
the simple reason that this is the first case which has
not been solved.  Our hope is to discover some 
structure which can be used to attack the general case.
It turns out that the $N=13$ case leads to relations reminiscent of
``higher order modular forms,'' which are described in the next section.
In Section~\ref{sec:prior} we describe prior work and then
in Section~\ref{sec:13} we apply those methods to the case~$N=13$.

In recent work, Conrey, Odgers, Snaith, and the first author~\cite{CFOS} 
have used some of the relations in this paper along with a new
generalization of Weil's lemma to complete the proof for~$N=13$.

\section{Higher order modular forms}

Our discussion here is imprecise and will only convey the general
flavor of this new subject.  For details see~\cite{CDO, DKMO}. 

We first introduce some slightly simpler notation.  If
$f|\gamma=f$ then we have
\begin{equation}
\gamma \equiv 1 \bmod \Omega_f
\end{equation}
where $\Omega_f$ is the right ideal in the group ring $\C[GL(2,\R)]$
which annihilates~$f$, the action of matrices on $f$ being extended
linearly.
We will write $\gamma \equiv 1$ instead of $\gamma \equiv 1 \bmod \Omega_f$
throughout this paper.
Thus, if $f$ is a cusp form for the group $\Gamma$, then the invariance 
properties of $f$ can be written as
$f|(1-\gamma)=0$ for all $\gamma\in \Gamma$, or equivalently,
$1-\gamma \equiv 0$.
This notation will make it easier to describe the properties of
higher order modular forms.

If $f$ is a \emph{second order} cusp form for the group $\Gamma$, 
then $f$ satisfies the
relation
\begin{equation}\label{eqn:secondorder}
(1-\gamma_1)(1-\gamma_2)\equiv 0
\end{equation}
for all
$\gamma_1, \gamma_2 \in \Gamma$.  Similarly, third order modular
forms satisfy 
\begin{equation}\label{eqn:thirdorder}
(1-\gamma_1)(1-\gamma_2)(1-\gamma_3)\equiv 0,
\end{equation}
and so on.  Roughly speaking, if $f$ is an $n$th order modular
form then $f|(1-\gamma)$ is an $(n-1)$st order modular form.
There are additional conditions involving the cusps and the parabolic
elements of $\Gamma$, but our goal here is just to introduce
the general idea.  Indeed, it is nontrivial to determine
the proper technical conditions, see~\cite{CDO, DKMO}.

In connection with our exploration of Conjecture~\ref{theconjecture}, 
a condition of form~(\ref{eqn:secondorder})
will arise where $\gamma_1$ and $\gamma_2$ come from \emph{different} groups.
This first appeared in the original work of Weil on the converse theorem
involving functional equations for twists.  Specifically, the
relation~(\ref{eqn:secondorder}) arose where
$\gamma_2$ was \emph{elliptic of infinite order}.  The following
lemma applies:

\begin{lemma}\label{thm:weil}  Suppose $f$ is holomorphic in $\H$ and
$\varepsilon\in GL_2(\R)^+$ is elliptic.  If \hbox{$f|_k\varepsilon=f$},
then either $\varepsilon$ has finite order, or $f$ is constant.
\end{lemma}

This is known as ``Weil's lemma''~\cite{W}.  See also the discussion in
Section~7.4 of Iwaniec's book~\cite{Iw}.
By the lemma, if $\gamma_2$ is elliptic of infinite order then
(\ref{eqn:secondorder}) implies that actually $1-\gamma_1\equiv 0$,
which is the conclusion Weil sought.  

Denote by $S_k(\Gamma_1,\Gamma_2)$
the set of analytic functions (with appropriate technical
conditions) satisfying~(\ref{eqn:secondorder})
for all $\gamma_1\in\Gamma_1$ and $\gamma_2\in\Gamma_2$.
The above lemma says that if $\Gamma_2$ contains an elliptic
element of infinite order then
$S_k(\Gamma_1,\Gamma_2) = S_k(\Gamma_1)$.  Note that the
analyticity of $f$ is necessary, and an analogue of Weil's
converse theorem for Maass form $L$-functions has not been
proven in classical language.

In Section~\ref{sec:13} we will see that our assumptions on the
Fricke involution and the Hecke operators lead to
condition~(\ref{eqn:secondorder})
with $\gamma_1\in\Gamma_0(13)$ and $\gamma_2$ in some other
discrete group.  We also obtain higher order
conditions~(\ref{eqn:thirdorder})
where each $\gamma_j$ comes from a different group.
This suggests the following question:

\begin{question}What conditions on $\Gamma_1$ and $\Gamma_2$ ensure that
$S_k(\Gamma_1,\Gamma_2)$ is finite dimensional?  What conditions
imply that $S_k(\Gamma_1,\Gamma_2) = S_k(\Gamma_1)$?
\end{question}

Part of the problem is determining the appropriate technical
conditions to incorporate into the definition of~$S_k(\Gamma_1,\Gamma_2)$.
Even when $\Gamma_1=\Gamma_2$ this is nontrivial. See~\cite{CDO, DKMO}.

\section{Manipulating the Hecke Operators}\label{sec:prior}

In \cite{CF} results were obtained for various $N$ up to $N=23$.
The idea is to manipulate the relations $T\equiv 1$, $H_N\equiv \pm 1$ and 
$T_n\equiv a_n$ to obtain $\gamma\equiv 1 $ for all 
$\gamma$ in a generating set for $\Gamma_0(N)$.   We will describe
the cases of $N=5,7,9,11$ from~\cite{CF}, and then the remainder of the paper
will concern the interesting relationships that arose in our exploration
of the case~$N=13$.

We have the following generating sets:
\begin{eqnarray}
\Gamma_0(N)&=&\left\langle T,\ W_N,\ \mm{2}{-1}{-N}{\frac{N+1}{2}} \right\rangle,
\phantom{XXXXX} N=5,7,9 , \cr
\Gamma_0(11)&=&\left\langle T,\ W_{11},\ \mm{2}{-1}{{-11}}{6} , 
\mm{3}{-1}{-11}{4}  \right\rangle, \cr
\Gamma_0(13)&=&\left\langle T,\ W_{13},\ \mm{2}{-1}{{-13}}{7} , 
\mm{-3}{-1}{13}{4} , \mm{3}{-1}{13}{-4} \right\rangle ,
\end{eqnarray}
where
\begin{equation}
T=\mm 11{}1
\ \ \ \ \ \ \ \ \ \ \ 
\mathrm{and}
\ \ \ \ \ \ \ \ \ \ \ 
W_N=\mm 1{}{N}1 .
\end{equation}

The generator $T$ is for free
because we have assumed a Fourier expansion.  The generator
$W_N$ now follows from the Fricke relation, because $W_N=H_N T H_N$.
So for these groups we have two of the generators.
Note that this uses the functional equation, but not the Euler product.

In the next section we repeat the calculations from~\cite{CF} in
the cases $N=5,7,9,11$, and in the following sections we treat
the case $N=13$.

\subsection{Levels 5, 7, 9, and 11}

For every $N$ we obtain a new generator from $T_2$.
This will resolve the cases $N=5$, $7$, and~$9$.

\begin{lemma}[Lemma 2 of \cite{CF}]  If $H_N \equiv \pm 1$ and
$T_2\equiv a_2$ then
$$
\mm 2{-1}{-N}{\frac{N+1}{2}} 
\equiv 1 .
$$
\end{lemma}

\begin{proof}
Note that
$$
H_N^{-1} T_2 H_N 
=  \mm 1{}{}2 + \mm 2{}{}1 + \mm 2{}{-N}1.
$$
Since $H_N^{-1} T_2 H_N \equiv a_2 H_N^{-1} H_N \equiv a_2\equiv T_2$, we have:
$$
\mm 1{}{}2 + \mm 2{}{}1 + \mm 2{}{-N}1 \equiv 
\mm 2{}{}1 + \mm 1{}{}2 + \mm 11{}2.
$$

Canceling common terms from both sides we are left with
$$
\mm 2{}{-N}1 \equiv \mm 11{}2.
$$
Right multiplying by $\mm 11{}2^{-1}$ we have
$$
M_2:=\mm 2{-1}{-N}{\frac{N+1}{2}} \equiv 1.
$$
\end{proof}
The lemma provides the final generator for 
$\Gamma_0(5)$, $\Gamma_0(7)$, and  $\Gamma_0(9)$.

To obtain the final generator for $\Gamma_0(11)$ we will combine
the Hecke operators $T_3$ and $T_4$ 
For $T_3$ we have
\begin{eqnarray}\label{eqn:T3manipulation}
0 &\equiv & \mathstrut H_{N}(T_3-a_3)H_{N} -(T_3-a_3)  \cr
&=&
- \mm 1{1}{}3 - \mm 1{2}{}3 + \mm 3{}{-N}1 + \mm 3{}{-2N}1 \cr
&\equiv& 
- \mm 1{1}{}3 - \mm 1{-1}{}3 + \mm 3{}{-N}1 + \mm 3{}{N}1 , 
\end{eqnarray}
where the second step used
\begin{equation}
\mm 1{-1}{}{1}\equiv 1
\ \ \ \ \ \ \ 
\text{ and }
\ \ \ \ \ \ \ 
\mm 1{}{N}{1}\equiv 1.
\end{equation}

We can combine the terms in pairs using
$$ 
\mm {1}{a}{}{p} -  \mm {p}{}{Nb}{1} = \( 1- \mm {p}{-a}{Nb}{\frac{-Nab+1}{p}} \) \mm {1}{a}{}{p}
$$
to get
\begin{equation}\label{eqn:T3rel}
\left(1-\mm 3{-1}{-11}4 \right)\beta(1/3)
+\left(1- \mm 31{11}4 \right)\beta(-1/3)
\equiv 0 ,
\end{equation}
where $\beta(x)=\displaystyle{\mm{1}{x}{}{1}}$.
We will combine this with a relation obtained from $T_4$ .

Since $T_4$ and $T_2$ are not
independent, there is more than one way to proceed.  The calculation
which seems most natural to us begins with
\begin{eqnarray}
0&\equiv &
H_{N}(T_{4} -a_4 )H_{N} - (T_{4} -a_4 ) \cr
&&\mathstrut - \left[H_{N}(T_{2} -a_2 )H_{N} -(T_{2} -a_2 )\right]\mm 2{}{}1 \cr
&&\mathstrut - \left[H_{N}(T_{2} -a_2 )H_{N}-(T_{2} -a_2 )\right]\mm 1{}{}2 \cr
&= &
- \mm {1}{1}{}{4} + \mm{4}{}{-3N}{1}
 - \mm {1}{3}{}{4} + \mm {4}{}{-N}{1} .
\end{eqnarray}
Combining terms as in the $T_3$ case gives
\begin{equation}\label{eqn:T4rel}
\left(1- \mm 4{-1}{-11}3 \right)\beta(1/4) 
+\left(1- \mm 41{11}3 \right)\beta(-1/4)
\equiv 0. 
\end{equation}
Combining~(\ref{eqn:T3rel}) and~(\ref{eqn:T4rel}) we obtain
\begin{eqnarray}
1- \mm 3{-1}{-11}4 
&\equiv&\mathstrut-\left(1- \mm 31{11}4 \right)\beta\left(-\frac 23 \right)\cr 
&=&\mathstrut\left(1- \mm 4{-1}{-11}3 \right) \mm 31{11}4 \beta\left(-\frac 23 \right)\cr 
&\equiv&\mathstrut -\left(1- \mm 41{11}3 \right)\beta\left(-\frac{2}{4} \right)
\mm 31{11}4 \beta\left(-\frac 23 \right)\cr 
&=&\mathstrut\left(1-\mm 3{-1}{-11}4 \right)\mm 41{11}3 \beta\left(
-\frac 24 \right)\mm  31{11}4 \beta\left(-\frac 23 \right) .
\end{eqnarray}
However,
$$
\mm 41{11}3
\beta\left(-\frac{2}{4} \right)
\mm 31{11}4 \beta\left(-\frac 23 \right)
= \mm 1{-2/3}{11/2}{-8/3} 
$$
is elliptic but not of finite order.
So by Lemma~\ref{thm:weil},
$$ 
\mm 3{-1}{-11}4 \equiv 1 .
$$
This is the final generator for~$\Gamma_0(11)$.

\section{Level $13$, mimic previous methods}\label{sec:13}

We will mimic the method used for $\Gamma_0(11)$ for 
$\Gamma_0(13)$, but things will not work out as nicely.  
What will arise is an expression of the form~(\ref{eqn:secondorder}) that
appears in the definition of second order modular form.

\subsection{The case of $T_3$. }\label{sec:T3}

From $ T_3$ we obtain the following expression, which is analogous
to~(\ref{eqn:T3rel}), 
\begin{equation}
\left(1-\mm 3{-1}{13}{-4} \right)\beta(1/3)
+\left(1- \mm 3{1}{-13}{-4} \right)\beta(-1/3)
\equiv 0 .
\end{equation}
We manipulate this similarly to the example for $\Gamma_0(11)$:
\begin{eqnarray}
1-\mm 3{1}{-13}{-4}
&\equiv&\mathstrut-\left(1- \mm 3{-1}{13}{-4} \right)\beta\left(\frac 23 \right)\cr 
&=&\mathstrut\left(1- \mm {-4}{1}{-13}3 \right) \mm 3{-1}{13}{-4} 
                    \beta\left(\frac 23 \right)\cr 
&\equiv&\mathstrut\left(1- H_{13} \mm {-4}{1}{-13}3 H_{13} \right) 
          H_{13} \mm 3{-1}{13}{-4} 
                    \beta\left(\frac 23 \right)\cr 
&=& \mathstrut\left(1- \mm 3{1}{-13}{-4}\right) 
          H_{13} \mm 3{-1}{13}{-4} 
                    \beta\left(\frac 23 \right) .
\end{eqnarray}
So,
\begin{equation}
\left(1- \mm 3{1}{-13}{-4}\right) \left ( 1- \varepsilon_1\right) \equiv 0
\end{equation}
where
\begin{equation}
\varepsilon_1=   H_{13} \mm 3{-1}{13}{-4} 
                    \beta\left(\frac 23 \right)
=  \mm {\sqrt{13}}{\frac{14}{3 \sqrt{13}}}{-3 \sqrt{13}}{- \sqrt{13}} .
\end{equation}
Since $\varepsilon_1$ is elliptic of order 2 
we cannot obtain anything from Lemma~\ref{thm:weil}.
However, we do have an expression of the form~(\ref{eqn:secondorder}) which
looks like the definition of a second order modular form.

\subsection{The case of $T_4$. }\label{sec:T4}
From $T_4$, again proceeding as in the $\Gamma_0(11)$ example,
we first have
\begin{equation}
\left(1- \mm 4{-1}{13}{-3} \right)\beta(1/4) 
+\left(1- \mm 41{-13}{-3} \right)\beta(-1/4)
\equiv 0. 
\end{equation}
Continuing exactly as above, this leads to
\begin{equation}
\left(1- \mm 3{1}{-13}{-4}\right) \left ( 1- \varepsilon_2\right) \equiv 0
\end{equation}
where
\begin{equation}
\varepsilon_2
= \mm {-\sqrt{13}}{\frac{-4}{ \sqrt{13}}}{\frac{7 \sqrt{13}}{2}}{\sqrt{13}} .
\end{equation}
Again $\varepsilon_2$ is elliptic of order 2.

\subsection{Combining $T_3$ and $T_4$. }
We can combine the two relationships to obtain 
\begin{equation}
0\equiv   \left[1- \mm{3}{1}{-13}{-4} \right] 
 \left (1 -\varepsilon \right)
\end{equation}
for any $\varepsilon$ in the group generated by $\varepsilon_1$ and $\varepsilon_2$,
and perhaps one of those elements will be elliptic of infinite order?
Unfortunately, this is not the case.  Note that
$$
\varepsilon_1 \varepsilon_2 = \mm{\frac{10}{3}}{\frac{2}{3}}{-\frac{13}{2}}{-1} ,
$$
which is hyperbolic.  Since $\varepsilon_1$ and $\varepsilon_2$ have order~2,
the group they generate contains only the elements
$(\varepsilon_1 \varepsilon_2)^n$ and $\varepsilon_2 (\varepsilon_1 \varepsilon_2)^n$,
so that group is discrete.

Although $T_3$ and $T_4$ were not sufficient to obtain the missing
generator, there are an infinite number of other Hecke operators to try.

\subsection{The case of $T_6$. } 

We now proceed with similar calculations with $T_6$.   We have
\begin{eqnarray}
0&\equiv& H_{13}(T_{6} -a_6 )H_{13} -(T_{6} -a_6 ) \cr
&&- \left[H_{13}(T_{2} -a_2 )H_{13} -(T_{2} -a_2 )\right] \mm 3{}{}1
- \left[H_{13}(T_{2} -a_2 )H_{13} -(T_{2} -a_2 )\right] \mm 1{}{}3 \cr
&&- \left[H_{13}(T_{3} -a_3 )H_{13} -(T_{3} -a_3 )\right]\mm 2{}{}1
- \left[H_{13}(T_{3} -a_3 )H_{13} -(T_{3} -a_3 )\right]\mm 1{}{}2 \cr
&\equiv &- \mm {1}{1}{}{6} + \mm{6}{}{-65}{1}
 - \mm {1}{5}{}{6} + \mm {6}{}{-13}{1} .
\end{eqnarray}
Using manipulations similar to those above gives
\begin{eqnarray}
0&\equiv& \left[-1+ \mm{6}{-1}{65}{11} \right] \mm{1}{1}{}{6}  + \left[- 1 + \mm{6}{-5}{-13}{11} \right] 
\mm{1}5{}6 \cr
&\equiv&\left[-1+ \mm{6}{-1}{65}{11} \right] \mm{1}{1}{}{6} , \cr
\end{eqnarray}
because $ \mm{6}{-5}{-13}{11} = M_2^{-1} H_{13}T^{-1} H_{13} T^{-1}$ so the 
second term on the first line is~$\equiv 0$.  So we have
$$
\mm{1}{1}{}{6}  + \mm{6}{}{-65}{1} \equiv 0 ,
$$
so
$$
\mm  6{-1}{-65}{11} \equiv 1
$$
This is not a new matrix because $\mm  6{-1}{-65}{11} =H_{13}T H_{13}T H_{13} M_2H_{13}$. 
That is, the above manipulations with $T_6$ produce results that can
be obtained from $T_2$.  

\subsection{Computer manipulation of Hecke operators}\label{sec:132}

The explicit manipulation of Hecke operators described in this paper are quite
tedious to do by hand, so we decided to make use of a computer.
We modified Mathematica to do calculations in the group ring
$\C[SL(2,\R)]$, made functions for the Hecke operators,
automated manipulations that occur repeatedly (such as the first step
in every example in the previous section of this paper), and 
implemented some crude simplifications procedures.

For the simplification procedures, we sought to automate the discovery,
for example, that if $T\equiv1$, $H_{13}\equiv \pm 1$, and $M_2\equiv 1$,
then
\begin{equation}
-1+ \mm{6}{-1}{65}{11} \equiv 0,
\end{equation}
as we saw at the end of the previous section.  
Our approach was to put all of the matrices in each expression
in ``simplest form'' by considering all products (on the left)
with, for example, fewer than 6 matrices where are known 
to be $\equiv 1$, and then keeping the representative which
has the smallest entries.  This idea worked surprisingly well.

We also implemented a ``factorization'' function which would do the
(trivial) calculation to check such things as whether 
$1-\gamma_1-\gamma_2+\gamma_3$ was of the form
$(1-\gamma_1)(1-\gamma_2)$ or $(1-\gamma_2)(1-\gamma_1)$.

\subsection{The case of $T_7$. }

Calculations with $T_7$ yield interesting results.   We have
\begin{eqnarray}
0 &\equiv & H_{13}(T_7 -a_7 )H_{13} -(T_7 -a_7 ) \cr
&\equiv& - \mm {1}{2}{}{7} + \mm{7}{}{-52}{1} - \mm {1}{3}{}{7} + \mm {7}{}{-65}{1} \cr
&&\mathstrut - \mm {1}{4}{}{7} + \mm {7}{}{-26}{1} - \mm {1}{5}{}{7} + \mm {7}{}{-39}{1} . \cr
\end{eqnarray}
Note that the expression on the right consists of 4~pair of matrices, as
opposed to the 6~pair that one would expect to obtain from~$T_7$.
This is because 
two pair 
canceled during simplification.

It turns out that the right side of the above expression factors as
\begin{eqnarray}
&& \left[- 1+ \mm{-3}{1}{-13}{4} \right] \mm 12{}7  
 +\left[- 1+ \mm{7}{4}{-65}{-37} \right] \mm 1{-4}{}7 \cr
&&\hskip 0.3cm  + \left[- 1+ \mm{7}{-4}{-26}{15} \right] \mm 1{4}{}7  
  +\left[- 1+ \mm{3}{1}{-13}{-4} \right]\mm 1{-2}{}7 \cr
&&\hskip 0.7cm= \left[ -\mm {3}{1}{-13}{-4} +1 \right] \mm {3}{1}{-13}{-4} ^{-1} H_{13} \mm {1}{2}{}{7} \cr
&&\hskip 1.2cm+ \left[ -\mm {3}{1}{-13}{-4} +1 \right] \mm {3}{1}{-13}{-4} ^{-1} H_{13} 
							\mm 7{-1}{-13}{2} ^{-1} \mm 1{-4}{}7 \cr
&&\hskip 1.2cm +\left[- 1+ \mm{3}{1}{-13}{-4} \right] \mm {7}{1}{13}{2} ^{-1} \mm 1{4}{}7  +
					 \left[- 1+ \mm{3}{1}{-13}{-4} \right]\mm 1{-2}{}7 \cr
&&\hskip 0.7cm=  \left[- 1+ \mm{3}{1}{-13}{-4} \right]\cr 
&&\hskip 1.2cm \left(-\mm{-\sqrt{13}}{\frac{2}{\sqrt{13}}}{3 \sqrt{13}}{-\sqrt{13}}
 - \mm {2\sqrt{13}}{\frac{1}{\sqrt{13}}}{-7\sqrt{13}}{} + \mm 21{-13}{-3} 
+\mm {1}{-2}{}{7} \right). \cr
\end{eqnarray}

We can right multiply by the inverse of any of the four matrices in the second
factor to rewrite this in the form $(1- \gamma) (1+ A-B -C)$.
For no good reason we choose the first term, giving
\begin{eqnarray}
0&\equiv& \left[ 1- \mm{3}{1}{-13}{-4} \right] \cr
&&\hskip 1cm \times \left( 1 +\mm {\frac{29}{7}}{\frac{5}{7}}{-13}{-2}
-\mm {\frac{5\sqrt{13}}{7}}{\frac{17}{7\sqrt{13}}}{\frac{-22\sqrt{13}}{7}}{\frac{-5\sqrt{13}}{7}}
-\mm {\frac{5\sqrt{13}}{7}}{\frac{24}{7\sqrt{13}}}{-3\sqrt{13}}{-\sqrt{13}} \right) \cr
&=& \left[ 1- \mm{3}{1}{-13}{-4} \right] \( 1+A-B-C\),
\end{eqnarray}
say.
This expression factors further.  Specifically, one can check that 
$A=C B$, so we have
\begin{equation}
0\equiv 
\left[ 1- \mm{3}{1}{-13}{-4} \right] (1-C)(1-B)
\end{equation}
Unfortunately, $B^2=1$,  so we cannot immediately
cancel the final factor to reduce to a second-order type expression.
It would be good if that happened, because we would have another matrix
to combine with the $\varepsilon_1$ and $\varepsilon_2$ from 
Sections~\ref{sec:T3} and~\ref{sec:T4}.

However, there is a curious benefit to having $B^2=1$, for we also
have $A B = C$, so
\begin{equation}
0\equiv 
\left[ 1- \mm{3}{1}{-13}{-4} \right] (1-A)(1-B).
\end{equation}

Note that if $B^2=1$, independent of any conditions on
$A$ and $C$, then
$( 1+A-B-C)(1+B)=(1-C A^{-1})(1+A B A^{-1}) A $, so
\begin{equation}\label{eqn:order2factorization}
0\equiv 
\left[ 1- \mm{3}{1}{-13}{-4} \right] (1-C A^{-1})(1+A B A^{-1}),
\end{equation}
which is almost a third-order condition.  Such expressions arise
whenever we have an order-2 matrix, so some types of factorization
are not a surprise.  In the particular case at hand, 
$C A^{-1}=A B A^{-1}$, which has order~2, so
$(1-C A^{-1})(1+A B A^{-1})=0$ and~(\ref{eqn:order2factorization})
contains absolutely no information.  Perhaps one should think that 
if $B^2=1$ then there always is some factorization, for
either~(\ref{eqn:order2factorization}) is 
nontrivial, or the expression factors nontrivially in another way.


\subsection{A few other cases}

From $T_{10}$ we get
\begin{eqnarray}
0&\equiv&\left[ 1- \mm{3}{1}{-13}{-4} \right] \cr
&&\hskip 1cm \times \left(1+
\mm {\frac{21}{5}}{\frac{2}{5}}{-13}{-1}
-\mm {\frac{2\,{\sqrt{13}}}{5}}{\frac{7}{5\,{\sqrt{13}}}}{\frac{-11\,{\sqrt{13}}}{5}}{\frac{-2\,{\sqrt{13}}}{5}}
-\mm {\frac{4\,{\sqrt{13}}}{5}}{\frac{19}{5\,{\sqrt{13}}}}{-3\,{\sqrt{13}}}{-{\sqrt{13}}} \right) \cr
&=&\left[ 1- \mm{3}{1}{-13}{-4} \right] (1+A-B-C),
\end{eqnarray}
say.  Again $A=CB$ and $B^2=1$, so we obtain two factorizations.

From $T_{15}$ we get
\begin{eqnarray}
0&\equiv&\left[ 1-  \mm{3}{1}{-13}{-4} \right] \cr 
&&\hskip 1cm  \times \left(1
+\mm {\frac{16}{5}}{1}{- \frac{117}{5} }{-7}
 -\mm {4\,{\sqrt{13}}}{\frac{15}{{\sqrt{13}}}}{\frac{-209\,{\sqrt{13}}}{15}}{-4\,{\sqrt{13}}}
-\mm {\frac{17\,{\sqrt{13}}}{15}}{\frac{4}{{\sqrt{13}}}}{\frac{-59\,{\sqrt{13}}}{15}}{-{\sqrt{13}  }} \right), \cr
\end{eqnarray}
which again factors in the same two ways.

From $T_{9}$ we get
\begin{eqnarray}
0&\equiv&\left[ 1- \mm{3}{1}{-13}{-4} \right] \cr
&&\hskip 1cm  \times \left(1 
+\mm {\frac{10}{3}}{1}{-\frac{13}{3}}{-1}
-\mm {2{\sqrt{13}}}{\frac{9}{{\sqrt{13}}}}{\frac{-53\,{\sqrt{13}}}{9}}{-2\,{\sqrt{13}\ }}
-\mm {\frac{7\,{\sqrt{13}}}{9}} {\frac{4}{{\sqrt{13}}}} {\frac{-25\,{\sqrt{13}}}{9}} {-{\sqrt{13}}} \right). \cr
\end{eqnarray}
which again factors in the same two ways.

It would be helpful to understand the underlying reason why
these expressions factor.

More time on the computer should produce more relations, but it is
not clear how they will combine to produce the desired result.
It would be interesting if the relations could build to the point where 
one could reduce higher order relations to lower order ones, 
which could then combine with previously found relations to 
cause additional cancellation, and so on, reducing down to 
the one missing generator for~$\Gamma_0(13)$.  It would be
more satisfying if one could find manipulations which produce any
specific matrix, as one does in the proof of Weil's converse
theorem.  

Our approach here is to look for factorizations 
$(1-\gamma)(1-\delta)(1-\varepsilon)\equiv 0$ in the
hopes of eliminating the last factor, perhaps because 
$\varepsilon$ is elliptic of infinite order.  In the case of expressions
that do not factor, it would be interesting to know if there are
cancellation laws beyond those implied by Weil's lemma.  That
is, are there conditions on $A$, $B$, $C$ such that
$f|(1+A-B-C) =0$ implies some apparently stronger condition on~$f$,
beyond those cases where $1+A-B-C$ factors and  Weil's lemma applies?

\subsection{A curiosity}

All the manipulations in this paper involve ``pairing up'' the terms
in a linear combination of matrices.  Usually there is a natural
way to do this, for one is hoping to produce matrices in $\Gamma_0(N)$.
However, it is possible to pair the matrices in different ways,
and one would like some justification for the choices and to know the
consequences of making the right (or wrong) choices.  This is
discussed extensively in~\cite{FW}.

We now give an example by repeating the analysis of Section~\ref{sec:prior}
making the wrong choices.
From~(\ref{eqn:T3manipulation}) with $N=11$ we have
\begin{equation}
\left( 1- \mm{3}{-1}{11}{-\frac{10}{3} } \right) \beta(1/3)
+ 
\left( 1- \mm{3}{1}{-11}{-\frac{10}{3} } \right) \beta(-1/3)
\equiv 0 ,
\end{equation}
where $\beta(x)=\displaystyle{\mm{1}{x}{}{1}}$.
Now doing manipulations exactly as in Section~\ref{sec:T3} we obtain
\begin{equation}
0\equiv \left( 1- \mm{3}{-1}{11}{-\frac{10}{3} } \right) (1-\varepsilon),
\end{equation}
where
\begin{equation}
\varepsilon = H_{11} \mm{3}{1}{-11}{-10/3} \beta(-2/3)= 
\mm {\sqrt{11}}{-\frac{4}{\sqrt{11}}}{3\sqrt{11}}{-\sqrt{11}} ,
\end{equation}
which has order~2. 

Note that the above manipulations cannot lead to
$\mm{3}{-1}{11}{-\frac{10}{3}} \equiv 1$.  
%
%
Indeed, 
if $p$ is prime, the group generated by $\Gamma_0(p)$ and $H_{p}$
is a maximal discrete subgroup of $SL(2,\R)$. So no manipulation
can lead to a new matrix which is~$\equiv 1$.
Yet, we do obtain additional
second order modular form type properties for newforms in
$S_k(\Gamma_0(11))$.   It is not clear what mechanism
will lead to the production of new matrices for~$N=13$, yet not
produce a contradiction when~$N=11$.

Using $T_4$ in the same way gives
\begin{equation}
0\equiv \left( 1- \mm{4}{-1}{11}{-\frac{5}{2} } \right) \left(1-
\mm {\sqrt{11}}{-\frac{3}{\sqrt{11}}}{4\sqrt{11}}{-\sqrt{11}}  \right),
\end{equation}
and from $T_6$ you get
\begin{equation}
0\equiv \left( 1- \mm{6}{-1}{11}{-\frac{5}{3} } \right) \left(1-
\mm {\sqrt{11}}{-\frac{2}{\sqrt{11}}}{6\sqrt{11}}{-\sqrt{11}}  \right),
\end{equation}
where the inner matrix is hyperbolic.

This illustrates that $f|(1-\varepsilon)(1-\delta)=0$ need not
imply $f$ is constant, and even having multiple independent relations of that
form is not sufficient.  In the case here, we have the above relations
in addition to $f|(1-\gamma)$ for all $\gamma\in\Gamma_0(11)$.  This
suggest that these ``second order'' conditions may be weaker than they
appear.

\end{document}